\documentclass[leqno,12pt]{amsart} 
\usepackage{amssymb} 
\usepackage{amsmath}
\usepackage{amsthm}
\usepackage{bm}
\usepackage{graphicx,psfrag,wrapfig}
\theoremstyle{plain} 
\newtheorem{theorem}{\indent\sc Theorem}[section] 
\newtheorem{lemma}[theorem]{\indent\sc Lemma}
\newtheorem{corollary}[theorem]{\indent\sc Corollary}

\theoremstyle{definition} 

\newtheorem{remark}[theorem]{\indent\sc Remark}

\newcommand{\affine}{\mathbb{C}}

\newcommand{\tangentvector}[1]{\partial / \partial #1}

\newcommand{\Dnorm}[3]{\left|\!\left| #1 \right|\!\right|_{\mathcal{C}^{#2} #3}}

\begin{document}

\title[On Holomorphic Curves in algebraic torus]
{On Holomorphic Curves in algebraic torus} %

\author[Masaki Tsukamoto]{Masaki Tsukamoto} %

\subjclass[2000]{ 
32H30}

\keywords{entire holomorphic curve in $(\mathbb{C}^*)^n$, polynomial growth}

\maketitle

\begin{abstract}
We study entire holomorphic curves in the algebraic torus, and show that 
they can be characterized by the ``growth rate'' of their derivatives.
\end{abstract}

\section{Introduction}
Let $z = x + y \sqrt{-1}$ be the natural coordinate in the complex plane $\affine$, and let $f(z)$ be an entire holomorphic function in the complex plane. Suppose that there are a non-negative 
integer $m$ and a positive constant $C$ such that 
\[ |f(z)| \leq C |z|^m, \quad (|z| \geq 1). \]
Then $f(z)$ becomes a polynomial with $\deg f(z) \leq m$.
This is a well-known fact in the complex analysis in one variable.
In this paper, we prove an analogous result for entire holomorphic curves in the algebraic torus $(\mathbb{C}^*)^n := (\affine \setminus \{0\})^n$.

Let $[z_0:z_1:\cdots :z_n]$ be the homogeneous coordinate in the complex projective space $\affine P^n$.
We define the complex manifold $X \subset \affine P^n$ by 
\[ X := \{ [1:z_1:\cdots :z_n] \in \affine P^n |\, z_i \neq 0, \, (1 \leq i \leq n) \} \cong (\mathbb{C}^*)^n .\]
$X$ is a natural projective embedding of $(\mathbb{C}^*)^n$.
We use the restriction of the Fubini-Study metric as the metric on $X$.

For a holomorphic map $f: \affine \to X$, we define its norm $|df|(z)$ by setting 
\begin{equation}\label{def: norm}
 |df|(z) := \sqrt{2}\, |df(\partial /\partial z)| \quad \text{for all $z \in \affine$} .
\end{equation}
Here $\tangentvector{z} = \frac{1}{2} \,(\tangentvector{x} - \sqrt{-1} \tangentvector{y})$, and 
the normalization factor $\sqrt{2}$
comes from $|\tangentvector{z}| = 1/\sqrt{2}$. 

The main result of this paper is the following.
\begin{theorem}\label{main theorem}
Let $f: \affine \to X$ be a holomorphic map. Suppose there are a non-negative integer $m$
and a positive constant $C$ such that 
\begin{equation}
|df|(z) \leq C |z|^m, \quad (|z|\geq 1).  \label{polynomial growth}
\end{equation}
Then there are polynomials $g_1(z)$, $g_2(z)$, $\cdots$, $g_n(z)$ 
with $\deg g_i(z) \leq m+1$, $(1\leq i \leq n)$, such that 
\begin{equation}
f(z) = [1: e^{g_1(z)} : e^{g_2(z)}: \cdots : e^{g_n(z)}].  \label{exp(m+1)}
\end{equation}

Conversely, if a holomorphic map $f(z)$ is expressed by $(\ref{exp(m+1)})$ with polynomials $g_i(z)$ of degree at most $m+1$, 
$f(z)$ satisfies the ``polynomial growth condition'' $(\ref{polynomial growth})$.
\end{theorem}
The direction (\ref{exp(m+1)}) $\Rightarrow$ (\ref{polynomial growth}) is easier, 
and the substantial part of the argument is the direction (\ref{polynomial growth}) $\Rightarrow$ 
(\ref{exp(m+1)}).

If we set $m=0$ in the above, we get the following corollary.
\begin{corollary}\label{BD}
Let $f:\affine \to X$ be a holomorphic map with bounded derivative, i.e., $|df|(z) \leq C$ for 
some positive constant $C$.  Then there are complex numbers $a_i$ and $b_i$, $(1\leq i \leq n)$, such that
\[ f(z) = [\, 1: e^{a_1 z + b_1}: e^{a_2 z + b_2}: \cdots : e^{a_n z + b_n}\, ] . \]
\end{corollary}
This is the theorem of [BD, Appendice]. The author also proves this in [T, Section 6].

\begin{remark}
The essential point of Theorem \ref{main theorem} is the statement that 
the degrees of the polynomials $g_i(z)$ are at most $m+1$.
Actually, it is easy to prove that if $f(z)$ satisfies the condition (\ref{polynomial growth})
then $f(z)$ can be expressed by (\ref{exp(m+1)}) with polynomials $g_i(z)$ of degree at most 
$2m+2$. (See Section 4.) 
\end{remark}
Theorem \ref{main theorem} states that holomorphic curves in $X$ can be characterized by the growth 
rate of their derivatives. We can formulate this fact more clearly as follows;

Let $g_1(z),\, g_2(z),\, \cdots,\, g_n(z)$ be polynomials, and define $f:\affine \to X$  by 
(\ref{exp(m+1)}).
We define the integer $m \geq -1$ by setting 
\begin{equation}\label{def: max deg}
 m+1 := \max (\deg g_1(z), \deg g_2(z), \cdots , \deg g_n(z) ) .
\end{equation}
We have $m = -1$ if and only if $f$ is a constant map. $m$ can be obtained as the growth rate of $|df|$:
\begin{theorem}\label{thm: order of |df|}
If $m\geq 0$, we have 
\[ \limsup_{r\to \infty} \frac{\max_{|z| =r} \log |df|(z)}{\log r} = m. \]
\end{theorem}
\begin{corollary}
Let $\lambda$ be a non-negative real number, and let $[\lambda ]$ be the maximum integer not greater than
$\lambda$. Let $f:\affine \to X$ be a holomorphic map, and suppose that there is a positive constant 
$C$ such that 
\begin{equation}\label{real growth rate}
 |df| (z) \leq C |z|^{\lambda} , \quad (|z| \geq 1). 
\end{equation}
Then we have a positive constant $C'$ such that 
\[ |df|(z) \leq C' |z|^{[\lambda ]}, \quad (|z| \geq 1). \]
\end{corollary}
\begin{proof}
If $f$ is a constant map, the statement is trivial. Hence we can suppose $f$ is not constant.
From Theorem \ref{main theorem}, $f$ can be expressed by (\ref{exp(m+1)}) with polynomials $g_i(z)$ of degree at most 
$[\lambda ] + 2$. Since $f$ satisfies (\ref{real growth rate}), we have 
\[ \limsup_{r\to \infty} \frac{\max_{|z| =r} \log |df|(z)}{\log r} \leq \lambda .\]
From Theorem \ref{thm: order of |df|}, this shows $\deg g_i(z) \leq [\lambda ] +1$ for all $g_i(z)$.
Then, Theorem \ref{main theorem} gives the conclusion.
\end{proof}


\section{Proof of (\ref{exp(m+1)}) $\Rightarrow$ (\ref{polynomial growth})}
Let $f:\affine \to X$ be a holomorphic map. From the definition of $X$, we have holomorphic maps
$f_i : \affine \to \affine ^*$, $(1\leq i \leq n)$,
 such that $f(z) = [1: f_1(z) : \cdots : f_n(z) ]$.
The norm $|df|(z)$ in (\ref{def: norm}) is given by 
\begin{equation}\label{def: Fubini-Study}
|df|^2 (z) = \frac{1}{4\pi} \Delta \log \left( 1+ \sum_{i=1}^n |f_i(z)|^2\right) , \quad 
(\Delta := \frac{\partial^2}{\partial x^2} + \frac{\partial^2}{\partial y^2}
= 4\, \frac{\partial ^2}{\partial z \partial \bar{z}} ).
\end{equation}
Suppose that $f$ is expressed by (\ref{exp(m+1)}), i.e., $f_i(z) = \exp (g_i(z))$ with a polynomial 
$g_i(z)$ of degree $\leq m+1$. We will repeatedly use the following calculation in this paper.
\begin{equation} \label{estimate of |df|^2 by f_i and f_i/f_j}
 \begin{split}
 |df|^2 &= \frac{1}{\pi}\left[ \frac{\sum_i |f'_i|^2}{(1 + \sum_i |f_i|^2)^2} 
+ \frac{\sum_{i<j} |g'_i - g'_j|^2 |f_i|^2 |f_j|^2}{(1 + \sum_i |f_i|^2)^2} \right] , \\  
 &\leq \frac{1}{\pi} \left[ \sum_i \frac{|f'_i|^2}{(1 + |f_i|^2)^2} 
 + \sum_{i<j} \frac{|g'_i - g'_j|^2 |f_i|^2 |f_j|^2}{(|f_i|^2 + |f_j|^2)^2}\right] , \\
 &= \frac{1}{\pi} \left[ \sum_i \frac{|f'_i|^2}{(1 + |f_i|^2)^2} 
 + \sum_{i<j} \frac{|(f_i/f_j)'|^2}{(1 + |f_i/f_j|^2)^2} \right] ,\\
 &= \sum_i |df_i|^2 + \sum_{i<j} |d(f_i/f_j)|^2  .
 \end{split}
\end{equation}
Here we set 
\[ |df_i| := \frac{1}{\sqrt{\pi}} \frac{|f'_i|}{1 + |f_i|^2} \quad \text{and}\quad
|d(f_i/f_j)| := \frac{1}{\sqrt{\pi}}\frac{|(f_i/f_j)'|}{1 + |f_i/f_j|^2}. \]
These are the norms of the differentials of the maps $f_i,\, f_i/f_j: \affine \to \affine P^1$.

We have $f_i(z) = \exp (g_i(z))$ and $f_i(z)/ f_j(z) = \exp (g_i(z) -g_j(z))$, and the degrees of 
the polynomials $g_i(z)$ and $g_i(z) - g_j(z)$ are at most $m+1$. Then, the next Lemma gives the desired conclusion:
\[ |df| (z) \leq C |z|^m, \quad (|z| \geq 1), \]
for some positive constant $C$.
\begin{lemma}\label{lemma: |dexp(m+1)|}
Let $g(z)$ be a polynomial of degree $\leq m+1$, and set $h(z) := e^{g(z)}$.
Then we have a positive constant $C$ such that 
\[ |dh|(z) = \frac{1}{\sqrt{\pi}}\frac{|h'(z)|}{1 + |h(z)|^2} \leq C |z|^m, \quad (|z| \geq 1). \]
\end{lemma}
\begin{proof}
We have 
\[ \sqrt{\pi}\, |dh| = \frac{|g'|}{|h| + |h|^{-1}} \leq |g'| \min (|h|, |h|^{-1}) \leq |g'| .\]
Since the degree of $g'(z)$  is at most $m$, we easily get the conclusion.
\end{proof}


\section{Preliminary estimates}\label{section: preliminary estimates}
In this section, $k$ is a fixed positive integer. 

The following is a standard fact in the Nevanlinna theory.
\begin{lemma}\label{lemma: easy estimate}
Let $g(z)$ be a polynomial of degree $k$, and set $h(z) = e^{g(z)}$. Then we have a positive 
constant $C$ such that 
\[ \int_1^r \frac{dt}{t} \int_{|z| \leq t} |dh|^2(z)\, dxdy \leq C r^k , \quad (r \geq 1). \]
\end{lemma}
\begin{proof}
Since $|dh|^2 = \frac{1}{4\pi} \Delta \log (1 + |h|^2)$, Jensen's formula gives 
\begin{equation*}
 \int_1^r \frac{dt}{t} \int_{|z| \leq t} |dh|^2\, dxdy
 = \frac{1}{4\pi} \int_{|z| = r} \log ( 1 + |h|^2)\, d\theta 
    -  \frac{1}{4\pi} \int_{|z| = 1} \log ( 1 + |h|^2)\, d\theta . 
\end{equation*}
Here $(r, \theta)$ is the polar coordinate in the complex plane.
We have 
\[ \log (1+|h|^2) \leq 2\, |\mathrm{Re}\, g(z)| + \log 2 \leq C r^k, \quad (r := |z| \geq 1). \]
Thus we get the conclusion.
\end{proof}
Let $I$ be a closed interval in $\mathbb{R}$ and let $u(x)$ be a real valued function defined on $I$.
We define its $\mathcal{C}^1$-norm $\Dnorm{u}{1}{(I)}$ by setting 
\[ \Dnorm{u}{1}{(I)} := \sup_{x\in I} |u(x)| + \sup_{x\in I} |u'(x)|. \]
For a Lebesgue measurable set $E$ in $\mathbb{R}$, we denote its Lebesgue measure by $|E|$.
\begin{lemma}\label{lemma: cosx}
There is a positive number $\varepsilon$ satisfying the following$:$
If a real valued function $u(x) \in \mathcal{C}^1 [0, \pi]$ satisfies
\[ \Dnorm{u(x) - \cos x}{1}{[0,\, \pi ]} \leq \varepsilon ,\]
then we have 
\[ |u^{-1} ([-t, t]) | \leq 4 t \quad \text{for any $t \in [0, \varepsilon]$}. \]
\end{lemma}
\begin{proof}
The proof is just an elementary calculus.
For any small number $\delta >0$, if we choose $\varepsilon$ sufficiently small, we have 
\[ u^{-1} ([-t, t]) \subset [\pi /2 -\delta , \pi /2 + \delta ]. \]
Let $x_1$ and $x_2$ be any two elements in $u^{-1} ([-t, t])$. From the mean value theorem, we have 
$y \in [\pi /2 -\delta , \pi /2 + \delta ]$ such that 
\[ u(x_1) - u(x_2) = u'(y)\, (x_1 - x_2) .\]
From $\sin (\pi /2) = 1$, we can suppose that 
\[ |u'(y)| \geq 1/2 .\]
Hence 
\[ |x_1 - x_2| \leq 2 \, |u(x_1) - u(x_2)| \leq 4 t .\]
Thus we get 
\[ |u^{-1} ([-t, t])| \leq 4t .\]
\end{proof}
Using a scale change of the coordinate, we get the following.
\begin{lemma}\label{lemma: coskx}
There is a positive number $\varepsilon$ satisfying the following$:$
If a real valued function $u(x) \in \mathcal{C}^1 [0, 2\pi]$ satisfies
\[ \Dnorm{u(x) - \cos kx}{1}{[0,\, 2\pi ]} \leq \varepsilon ,\]
then we have 
\[ |u^{-1} ([-t, t]) | \leq 8 t \quad \text{for any $t \in [0, \varepsilon]$}. \]
\end{lemma}
\begin{proof}
\[u^{-1} ([-t, t]) = \bigcup_{j=0}^{2k-1} u^{-1} ([-t, t])\cap [j\pi /k,\, (j+1)\pi /k]. \] 
Applying Lemma \ref{lemma: cosx} to $u(x/k)$, we have 
\[ |u^{-1} ([-t, t])\cap [0, \pi /k]\, | \leq 4t/k.  \] 
In a similar way,
\[ |u^{-1} ([-t, t])\cap [j\pi /k,\, (j+1)\pi /k]\, | \leq 4t/k , \quad (j = 0, 1, \cdots , 2k-1). \]
Thus we get the conclusion.
\end{proof}
Let $E$ be a subset of $\affine$. For a positive number $r$, we set 
\begin{equation*}
 E(r) := \{ \theta \in \mathbb{R}/2\pi \mathbb{Z} |\, r e^{i\theta} \in E\} .
\end{equation*}
In the rest of the section, we always assume $k\geq 2$.
\begin{lemma}\label{lemma: monic}
Let $C$ be a positive constant, and let $g(z) = z^k + a_1 z^{k-1} + \cdots + a_k$ be a monic 
polynomial of degree $k$. Set 
\[ E := \{ z \in \affine |\, |\mathrm{Re}\, g(z) | \leq C |z| \} . \]
Then we have a positive number $r_0$ such that 
\[ |E(r)| \leq 8C/ r^{k-1}, \quad (r \geq r_0). \] 
\end{lemma}
\begin{proof}
Set $v(z) := \mathrm{Re}\, (a_1 z^{k-1} + a_2 z^{k-2} + \cdots + a_k)$. 
Then we have 
\[ |\mathrm{Re}\, g(re^{i\theta}) | \leq Cr \quad \Longleftrightarrow 
\quad |\cos k\theta + v(re^{i\theta})/r^k| \leq C/r^{k-1} .\]
Set $u(\theta) := \cos k\theta + v(re^{i\theta})/r^k$.
It is easy to see that 
\[ \Dnorm{u(\theta ) - \cos k\theta }{1}{[0, 2\pi]} \leq \mathrm{const}/r , \quad (r \geq 1). \]
Then we can apply Lemma \ref{lemma: coskx} to this $u(\theta)$, and we get 
\[ |E(r)| = |u^{-1}([-C/r^{k-1},\, C/ r^{k-1}]) | \leq 8C /r^{k-1} , \quad (r \gg 1).\]
\end{proof}
The following is the key lemma.
\begin{lemma}\label{lemma: key}
Let $g(z) = a_0 z^k + a_1 z^{k-1} + \cdots + a_k$ be a polynomial of degree $k$, $(a_0 \neq 0)$.
Set 
\[ E := \{ z\in \affine |\, |\mathrm{Re}\, g(z)| \leq |z| \} .\]
Then we have a positive number $r_0$ such that 
\[ |E(r)| \leq \frac{8}{|a_0|r^{k-1}} ,\quad (r \geq r_0) . \]
\end{lemma}
\begin{proof}
Let $\arg a_0$ be the argument of $a_0$, and set $\alpha := \arg a_0 /k$.
We define the monic polynomial $g_1(z)$ by 
\[ g_1(z) := \frac{1}{|a_0|} g(e^{-i\alpha}z) = z^k + \cdots .\]
Then we have
\[ |\mathrm{Re}\, g(re^{i\theta})| \leq r \Longleftrightarrow 
|\mathrm{Re}\, g_1(re^{i(\theta + \alpha )})| \leq r/|a_0| .\]
Hence the conclusion follows from Lemma \ref{lemma: monic}.
\end{proof}
\begin{lemma}\label{lemma: integration over E^c}
Let $g(z)$ be a polynomial of degree $k$, and we define $E$ as in Lemma \ref{lemma: key}.
Set $h(z) := e^{g(z)}$. Then we have 
\[ \int_{\affine \setminus E} |dh|^2(z) \,dxdy < \infty .\] 
\end{lemma}
\begin{proof}
Since $|h| = e^{\mathrm{Re}\, g}$, the argument in the proof of Lemma \ref{lemma: |dexp(m+1)|} gives
\[ \sqrt{\pi}\, |dh| \leq |g'|\min (|h|, |h|^{-1}) \leq |g'|\, e^{-|\mathrm{Re}\, g|}. \]
$g'(z)$ is a polynomial of degree $k-1$, and we have $|\mathrm{Re}\, g| > |z|$ for 
$z \in \affine \setminus E$. 
Hence we have a positive constant $C$ such that 
\[ |dh|(z) \leq C |z|^{k-1} e^{-|z|}, \quad 
\text{if $z \in \affine \setminus E$ and $|z|\geq 1$}. \]
The conclusion follows from this estimate.
\end{proof}


\section{Proof of  (\ref{polynomial growth}) $\Rightarrow$ (\ref{exp(m+1)})}\label{section: proof of main theorem}
Let $f = [1:f_1:f_2:\cdots :f_n] : \affine \to X$ be a holomorphic map with 
$|df|(z)\leq C|z|^m$, $(|z|\geq 1)$. Since $\exp : \affine \to \affine ^*$ is the universal covering, we have 
entire holomorphic functions $g_i(z)$ such that $f_i(z) = e^{g_i(z)}$.
We will prove that all $g_i(z)$ are polynomials of degree $\leq m+1$. 
The proof falls into two steps. In the first step, we prove all $g_i(z)$ are polynomials.
In the second step, we show $\deg g_i(z) \leq m+1$. The second step is the harder part of the 
proof. 

Schwarz's formula gives\footnote{The idea of using Schwarz's formula is due to [BD, Appendice]. The author 
gives a different approach in [T, Section 6].} 
\[ \pi r^k g_i^{(k)}(0) = k! \int_{|z|=r} \mathrm{Re}\, (g_i(z))\, e^{-k\sqrt{-1}\theta} d\theta 
= k! \int_{|z| = r} \log |f_i(z)|\, e^{-k\sqrt{-1}\theta} d\theta , \quad (k\geq 1).\]
We have 
\[ \left| \log |f_i|\right| \leq \log (|f_i| + |f_i|^{-1}) = \log (1 + |f_i|^2) - \log |f_i|
\leq \log (1+ \sum |f_j|^2) - \log |f_i|.  \]
Hence 
\[ \pi r^k |g_i^{(k)}(0)| \leq k! \int_{|z|=r} \log (1+ \sum |f_j|^2)\, d\theta 
- k! \int_{|z|=r}\log |f_i|\, d\theta . \]
Since $\log |f_i| = \mathrm{Re}\, g_i(z)$ is a harmonic function, the second term in the above is equal to 
the constant $-2\pi k!\, \mathrm{Re}\, g_i(0)$. 
Since $|df|^2 = \frac{1}{4\pi}\Delta \log (1 + \sum |f_j|^2)$, Jensen's formula gives 
\[\frac{1}{4\pi} \int_{|z|=r} \log (1+ \sum |f_j|^2)\, d\theta - 
\frac{1}{4\pi} \int_{|z|=1} \log (1+ \sum |f_j|^2)\, d\theta
= \int_1^r \frac{dt}{t} \int_{|z|\leq t} |df|^2(z) \, dxdy .\]
Thus we get 
\begin{equation}\label{estimate of g_i^{(k)}(0)}
\frac{r^k}{4k!}|g_i^{(k)}(0)| \leq \int_1^r \frac{dt}{t} \int_{|z|\leq t} |df|^2(z) \, dxdy + 
\mathrm{const}.
\end{equation}
Since $|df|(z) \leq C |z|^m$, $(|z|\geq 1)$, this shows $g_i^{(k)}(0) = 0$ for $k\geq 2m + 3$. Hence $g_i(z)$ are polynomials.

Next we will prove $\deg g_i(z) \leq m+1$.
We define $E_i, \, E_{ij} \subset \affine$, $(1\leq i \leq n, 1\leq i<j\leq n)$, by setting 
\begin{align*}
&\deg g_i(z) \leq m+1 \Longrightarrow E_i := \emptyset , \\
&\deg g_i(z) \geq m+2 \Longrightarrow E_i := \{ z\in \affine | \, |\mathrm{Re}\, g_i(z)| \leq |z| \} ,\\
&\deg (g_i(z) -g_j(z)) \leq m+1 \Longrightarrow E_{ij} := \emptyset ,\\
&\deg (g_i(z) -g_j(z)) \geq m+2 \Longrightarrow 
E_{ij} := \{ z\in \affine | \, |\mathrm{Re}\, (g_i(z) - g_j(z)) | \leq |z| \}.
\end{align*}
We set $E := \bigcup_{i} E_i \cup \bigcup_{i<j} E_{ij}$. 
Then we have $E(r) =\bigcup_{i} E_i(r) \cup \bigcup_{i<j} E_{ij}(r)$ for $r>0$.
From Lemma \ref{lemma: key}, we have positive constants $r_0$ and $C'$ such that 
\begin{equation}
|E(r)| \leq C'/r^{m+1}, \quad (r\geq r_0). \label{estimate of E(r)}
\end{equation}
We have 
\begin{equation}\label{E and E^c}
 \begin{split}
  \int_1^r \frac{dt}{t} \int_{|z|\leq t}& |df|^2(z) \, dxdy \\
  &= \int_1^r \frac{dt}{t} \int_{E\cap \{ |z|\leq t\}} |df|^2(z) \, dxdy
  + \int_1^r \frac{dt}{t} \int_{E^c \cap \{ |z|\leq t\}} |df|^2(z) \, dxdy .
 \end{split}
\end{equation}
Using (\ref{estimate of E(r)}) and $|df|(z)\leq C |z|^m$, $(|z|\geq 1)$, 
we can estimate the first term in (\ref{E and E^c}) as follows:
\[ \int_{E\cap \{ 1\leq |z|\leq t\}} |df|^2(z) \, dxdy 
\leq C^2 \int_{E\cap \{ 1\leq |z|\leq t\}} r^{2m+1} \, dr d\theta = C^2 \int_1^t r^{2m+1}|E(r)| dr .\]
If $t \geq r_0$, we have 
\[ \int_{r_0}^t r^{2m+1}|E(r)| dr \leq C' \int_{r_0}^t r^m dr = \frac{C'}{m+1} t^{m+1} -  \frac{C'}{m+1} {r_0}^{m+1}.\]
Thus 
\begin{equation}\label{estimate of the first term}
\int_1^r \frac{dt}{t} \int_{E\cap \{ |z|\leq t\}} |df|^2(z) \, dxdy 
\leq \mathrm{const} \cdot r^{m+1}, \quad (r \geq 1).
\end{equation}

Next we will estimate the second term in (\ref{E and E^c}) by using the inequality (\ref{estimate of |df|^2 by f_i and f_i/f_j}) given in Section 2:
\[ |df|^2 \leq \sum_i |df_i|^2 + \sum_{i<j} |d(f_i/f_j)|^2 .\]
If $\deg g_i(z) \leq m+1$, Lemma \ref{lemma: easy estimate} gives 
\[ \int_1^r \frac{dt}{t} \int_{E^c \cap \{ |z|\leq t\}} |df_i|^2(z) \, dxdy
\leq \int_1^r \frac{dt}{t} \int_{|z|\leq t} |df_i|^2(z) \, dxdy
\leq \mathrm{const} \cdot r^{m+1} .\]
If $\deg g_i(z) \geq m+2$, Lemma \ref{lemma: integration over E^c} gives 
\[ \int_{E^c \cap \{ |z|\leq t\}} |df_i|^2(z) \, dxdy 
\leq \int_{E_i^c \cap \{ |z|\leq t\}} |df_i|^2(z) \, dxdy \leq \mathrm{const} .\]
The terms for $|d(f_i/f_j)|$ can be also estimated in the same way, and we get 
\begin{equation}\label{estimate of the second term}
 \int_1^r \frac{dt}{t} \int_{E^c \cap \{ |z|\leq t\}} |df|^2(z) \, dxdy 
\leq \mathrm{const}\cdot r^{m+1} , \quad (r\geq 1).
\end{equation}
From (\ref{E and E^c}), (\ref{estimate of the first term}), (\ref{estimate of the second term}), we get 
\[  \int_1^r \frac{dt}{t} \int_{|z|\leq t} |df|^2(z) \, dxdy \leq \mathrm{const}\cdot r^{m+1},
\quad (r \geq 1) .\]
From (\ref{estimate of g_i^{(k)}(0)}), this shows $g_i^{(k)}(0) = 0$ for $k\geq m+2$. 
Thus $g_i(z)$ are polynomials with $\deg g_i(z) \leq m+1$.
This concludes the proof of Theorem \ref{main theorem}.

\section{Proof of Theorem \ref{thm: order of |df|} and a corollary}
\subsection{Proof of Theorem \ref{thm: order of |df|}}
The proof of Theorem \ref{thm: order of |df|} needs the following lemma.
\begin{lemma}\label{lem: preparation}
Let $k\geq 1$ be an integer, and let $\delta$ be a real number satisfying $0< \delta <1$.
Let $g(z) = a_0 z^k + a_1 z^{k-1} + \cdots + a_k$ be a polynomial of degree $k$, $(a_0\neq 0)$.
We set $h(z) := e^{g(z)}$ and define $E \subset \affine$ by 
\[ E := \{ z\in \affine |\, |\mathrm{Re}\, g(z)| \leq |z|^{\delta}\}.\]
Then we have 
\[ \int_{\affine \setminus E} |dh|^2 < \infty ,\]
and there is a positive number $r_0$ such that 
\[ |E(r)| \leq \frac{8}{|a_0|r^{k-\delta}}, \quad (r\geq r_0) .\]
\end{lemma}
\begin{proof}
This can be proven by the methods in Section \ref{section: preliminary estimates}.
We omit the detail.
\end{proof}
Let $g_1(z), \, g_2(z), \, \cdots ,\, g_n(z)$ be polynomials, and define the holomorphic map $f: \affine \to X$ and  the integer $m\geq -1$ by (\ref{exp(m+1)}) and (\ref{def: max deg}).
Here we suppose $m\geq 0$, i.e., $f$ is not a constant map. 
We will prove Theorem \ref{thm: order of |df|}. 

From Theorem \ref{main theorem}, we have 
\[ |df|(z) \leq \mathrm{const}\cdot |z|^m , \quad (|z| \geq 1).\]
It follows 
\[ \limsup_{r\to \infty} \frac{\max_{|z| =r} \log |df|(z)}{\log r} \leq m .\]
We want to prove that this is actually an equality. Suppose 
\[ \limsup_{r\to \infty} \frac{\max_{|z| =r} \log |df|(z)}{\log r} \lvertneqq m .\]
Then, if we take $\varepsilon >0$ sufficiently small, we have a positive number $r_0$ such that 
\begin{equation}\label{bound of contradiction}
|df|(z) \leq |z|^{m-\varepsilon} ,\quad (|z| \geq r_0).
\end{equation}
Schwarz's formula gives the inequality (\ref{estimate of g_i^{(k)}(0)}):
\begin{equation}\label{estimate by Schwarz}
\frac{r^k}{4k!}|g_i^{(k)}(0)| \leq \int_1^r \frac{dt}{t} \int_{|z|\leq t} |df|^2(z) \, dxdy + 
\mathrm{const}, \quad (k\geq 0).
\end{equation}
Let $\delta$ be a positive number such that $0 < \delta < 2\varepsilon$. We define $E_i$ and $E_{ij}$, 
$(1\leq i \leq n, \, 1\leq i<j\leq n)$, by setting
\begin{align*}
&\deg g_i(z) \leq m \Longrightarrow E_i := \emptyset , \\
&\deg g_i(z) = m+1 \Longrightarrow E_i := \{ z\in \affine | \, |\mathrm{Re}\, g_i(z)| \leq |z|^{\delta} \} ,\\
&\deg (g_i(z) -g_j(z)) \leq m \Longrightarrow E_{ij} := \emptyset ,\\
&\deg (g_i(z) -g_j(z)) = m+1 \Longrightarrow 
E_{ij} := \{ z\in \affine | \, |\mathrm{Re}\, (g_i(z) - g_j(z)) | \leq |z|^{\delta} \}.
\end{align*}
We set $E := \bigcup_{i} E_i \cup \bigcup_{i<j} E_{ij}$. Then, if we take $r_0$ sufficiently large, we have 
\begin{equation}\label{estimate of E(r) ver.2}
 |E(r)| \leq \mathrm{const}/r^{m+1-\delta} , \quad (r\geq r_0) .
 \end{equation}
We have 
\begin{equation*}
 \begin{split}
  \int_1^r \frac{dt}{t} \int_{|z|\leq t}& |df|^2(z) \, dxdy \\
  &= \int_1^r \frac{dt}{t} \int_{E\cap \{ |z|\leq t\}} |df|^2(z) \, dxdy
  + \int_1^r \frac{dt}{t} \int_{E^c \cap \{ |z|\leq t\}} |df|^2(z) \, dxdy .
 \end{split}
\end{equation*}
From (\ref{bound of contradiction}) and (\ref{estimate of E(r) ver.2}), the first term can be estimated as in Section \ref{section: proof of main theorem}:
\[ \int_1^r \frac{dt}{t} \int_{E\cap \{ |z|\leq t\}} |df|^2(z) \, dxdy 
\leq \mathrm{const}\cdot r^{m+1- (2\varepsilon -\delta )} , \quad (r\geq 1).\]
Using Lemma \ref{lem: preparation} and the inequality $|df|^2 \leq \sum_i |df_i|^2 + \sum_{i<j} |d(f_i/f_j)|^2 $, we can estimate the 
second term:
\[ \int_1^r \frac{dt}{t} \int_{E^c \cap \{ |z|\leq t\}} |df|^2(z) \, dxdy 
\leq \mathrm{const}\cdot \log r + \mathrm{const} \cdot r^m , \quad (r\geq 1). \]
Thus we get 
\[ \int_1^r \frac{dt}{t} \int_{E\cap \{ |z|\leq t\}} |df|^2(z) \, dxdy  
\leq \mathrm{const}\cdot r^{m+1- (2\varepsilon -\delta )} , \quad (r\geq 1).\]
Note that $2\varepsilon -\delta$ is a positive number. Using this estimate in (\ref{estimate by Schwarz}), we get 
\[ g_i^{(k)}(0) = 0, \quad (k\geq m+1). \]
This shows $\deg g_i(z) \leq m$. This contradicts the definition of $m$. 
\begin{remark}
The following is also true:
\[  \limsup_{r\to \infty} \frac{\max_{|z| \leq r} \log |df|(z)}{\log r} = m .\]
\end{remark}
\begin{proof}
We have 
\[ m= \limsup_{r\to \infty} \frac{\max_{|z| = r} \log |df|(z)}{\log r}
\leq \limsup_{r\to \infty} \frac{\max_{|z| \leq r} \log |df|(z)}{\log r}. \]
And we have $|df|(z) \leq \mathrm{const}\cdot |z|^m$, $(|z| \geq 1)$.
Thus 
\[\limsup_{r\to \infty} \frac{\max_{|z| \leq r} \log |df|(z)}{\log r} \leq m .\]
\end{proof}

\subsection{Order of the Shimizu-Ahlfors characteristic function}
For a holomorphic map $f: \affine \to X$, we define 
the Shimizu-Ahlfors characteristic function $T(r, f)$ by 
\[ T(r, f) := \int_1^r \frac{dt}{t} \int_{|z| \leq t} |df|^2(z) \, dxdy, \quad (r \geq 1). \]
The order $\rho_f$ of $T(r, f)$ is defined by 
\[ \rho_f := \limsup_{r\to \infty} \frac{\log T(r,f)}{\log r}. \]
$\rho_f$ can be obtained as the growth rate of $|df|$:
\begin{corollary}\label{cor: order and growth rate of |df|}
For a holomorphic map $f:\affine \to X$, we have 
\[ \rho_f < \infty \Longleftrightarrow 
\limsup_{r\to \infty} \frac{\max_{|z| =r} \log |df|(z)}{\log r} < \infty .\]
If these values are finite and $f$ is not a constant map, then we have 
\[ \rho_f = \limsup_{r\to \infty} \frac{\max_{|z| =r} \log |df|(z)}{\log r} +1 . \]
\end{corollary}
\begin{proof}
If $\rho_f < \infty$, the estimate (\ref{estimate by Schwarz}) shows that $f$ can be expressed by 
(\ref{exp(m+1)}) with polynomials $g_1(z), \, \cdots,\, g_n(z)$. Then we have
\[ \limsup_{r\to \infty} \frac{\max_{|z| =r} \log |df|(z)}{\log r} <\infty .\]
The proof of the converse is trivial.

Suppose $\rho_f < \infty$. Then we can express $f$ by $f(z) = [1:e^{g_1(z)}:\cdots :e^{g_n(z)}]$ with 
polynomials $g_1(z), \, \cdots,\, g_n(z)$. We set $f_i(z) := e^{g_i(z)}$, and define the integer $m$ by 
(\ref{def: max deg}).
Theorem \ref{thm: order of |df|} gives 
\[ \limsup_{r\to \infty} \frac{\max_{|z| =r} \log |df|(z)}{\log r} +1 = m+1. \]
The estimate (\ref{estimate by Schwarz}) gives 
\[ m+1 \leq \rho_f . \]
Since $|df| = \frac{1}{4\pi} \Delta \log (1 + \sum |f_i|^2)$, Jensen's formula gives 
\begin{equation*}
T(r, f)
 = \frac{1}{4\pi} \int_{|z| = r} \log ( 1 + \sum_i |f_i|^2)\, d\theta 
    -  \frac{1}{4\pi} \int_{|z| = 1} \log ( 1 + \sum_i |f_i|^2)\, d\theta .
\end{equation*}
Since $\deg g_i(z) \leq m+1$, we have 
\[ \log (1+ \sum_i |f_i|^2) \leq \mathrm{const} \cdot r^{m+1}, \quad (r\geq 1). \]
Hence 
\[ \rho_f \leq m+1 .\]
Thus we get 
\[ \rho_f = m+1 = \limsup_{r\to \infty} \frac{\max_{|z| =r} \log |df|(z)}{\log r} +1 .\]
\end{proof}
\begin{remark}
Of course, the statement of Corollary \ref{cor: order and growth rate of |df|} is not true for general 
entire holomorphic curves in the complex projective space $\affine P^n$.
For example, let $f:\affine \to \affine P^1$ be a non-constant elliptic function. Since $|df|$ is 
bounded all over the complex plane, we have 
\[ \limsup_{r\to \infty} \frac{\max_{|z| =r} \log |df|(z)}{\log r} = 0.\]
And it is easy to see 
\[ \rho_f = 2 \neq \limsup_{r\to \infty} \frac{\max_{|z| =r} \log |df|(z)}{\log r} + 1. \]
\end{remark}


\vspace{10mm}

\address{ Masaki Tsukamoto \endgraf
Department of Mathematics, Faculty of Science \endgraf
Kyoto University \endgraf
Kyoto 606-8502 \endgraf
Japan
}

\textit{E-mail address}: \texttt{tukamoto@math.kyoto-u.ac.jp}

\end{document}